\documentclass[12pt]{amsart}
\usepackage{latexsym}
\usepackage{amsmath,amsfonts,amssymb,epsfig}
\usepackage{amsthm,amscd}
\usepackage{mathrsfs} 
\usepackage{hyperref}
\usepackage{cancel}
\usepackage{wrapfig}

\numberwithin{equation}{section} 

\newtheorem{theorem}{Theorem}[section]
\newtheorem{definition}[theorem]{Definition}

\newtheorem{proposition}[theorem]{Proposition}

\newenvironment{klemma}[1][Key Lemma]{\begin{trivlist}
\item[\hskip \labelsep {\bfseries #1}]}{\end{trivlist}}

\newtheorem{rem}[theorem]{Remark}
\newenvironment{remark}{\begin{rem} \rm}{\end{rem}}
\newtheorem{ide}[theorem]{Idea}

\newtheorem{exa}[theorem]{Example}

\newtheorem{ques}[theorem]{Question}

\newtheorem{pro}[theorem]{Problem}

\newtheorem{spe}[theorem]{Speculation}

\font\bbb=msbm10 scaled 1100

\newcommand{\vol}{\nu}
\newcommand{\R}{\mbox{\bbb R}}       

\newcommand{\Z}{\mbox{\bbb Z}}
\newcommand{\Cnf}{\text{\rm Conf}}

\newcommand{\no}{\noindent}

\title[The third order helicity of magnetic fields via link maps II.]{The third order helicity of magnetic fields\\ via link maps II.} 

\author[R. Komendarczyk]{R. Komendarczyk$^\ast$}

\address{Department of Mathematics,
Tulane University,
6823 St. Charles Ave.,
New Orleans, LA 70118} \email{rako@tulane.edu}

\thanks{$^\ast$The author acknowledges the support  of DARPA \#FA9550-08-1-0386.}


\date{\today}


\begin{document}

\begin{abstract}
In this sequel we extend the derivation of the third order helicity to magnetic fields supported on unlinked domains in 3-space. The formula is expressed in terms of generators of the deRham cohomology of the configuration space of three points in $\R^3$, which is a more practical domain from the perspective of applications. It also admits an ergodic interpretation as an average asymptotic Milnor $\bar{\mu}_{123}$-invariant and allows us to obtain
the $L^2$-energy bound for the magnetic field. As an intermediate step we derive an integral formula for Milnor $\bar{\mu}_{123}$-invariant for parametrized Borromean links in $\R^3$.
\end{abstract} 

\maketitle

\section{Introduction}\label{sec:1}
In the recent work \cite{Kom-unlinked09} the author derived a new formula for the third order helicity 
$\mathsf{H}_{123}(B;\mathcal{T})$ of a volume preserving vector field $B$ supported on invariant unlinked
domains $\cup_i\mathcal{T}_i$ in the $3$-sphere $S^3$.  Here by {\em unlinked} we understand  disjoint compact handlebodies with smooth boundary such that every pair of 1-cycles  in 
$H_1(\mathcal{T}_i)$ and $H_1(\mathcal{T}_j)$, $i\neq j$ has a linking number zero, Figures \ref{fig:borro} and \ref{fig:invariant-sets} show examples of such domains.  Note that it is a much weaker property as unlinked in the standard sense of the word (see e.g. \cite{Prasolov-Sossinsky97}).

A purpose of this sequel is to derive a formula for $\mathsf{H}_{123}(B;\mathcal{T})$ for domains $\cup_i \mathcal{T}_i$ of the 3-space $\R^3$, Theorem \ref{th:main}$(i)$, as this setting is more natural from perspective of applications to fluid dynamics \cite{Priest84}.  The main theorems can be considered as an extension of Laurence and Stredulinsky results from \cite{Laurence-Stredulinsky00b, Laurence-Stredulinsky00a}
to vector fields supported on invariant unlinked handlebodies in $\R^3$. It may seem at first like a minor improvement since $S^3=\R^3\cup\{\infty\}$, and one could simply ``pull-back'' the formula obtained  in \cite{Kom-unlinked09} to $\R^3$. However, the new formula obtained here is qualitatively different, it involves familiar \emph{Green forms} $\{\omega_{1,2}, \omega_{2,3}, \omega_{3,1}\}$ representing generators of the cohomology ring of the configuration space $\Cnf_3(\R^3)$ of three points in $\R^3$. It also allows us to derive an  $L^2$-energy bound for $B$ (Theorem \ref{th:main2}) which involves the flat geometry rather than the spherical geometry, as in \cite{Kom-unlinked09}. In our Key Lemma, we obtain an integral for the $\bar{\mu}_{123}$-invariant of 3-component \emph{Borromean links} in $\R^3$, i.e. links with vanishing pairwise linking 
numbers. Note that the Borromean links are also known as \emph{homotopy Brunnian} \cite{Koschorke91, Koschorke97}.

Helicity invariants measure topological complexity of the flow and are relevant in the context of e.g. plasma physics where $B$ is a magnetic field frozen in the velocity field of plasma \cite{Priest84, Moffatt95,Moffatt85}. Most known helicity invariants, such as Woltjer's helicity \cite{Arnold86} or higher helicities introduced in \cite{Berger90, Hornig04} are vector field analogs of Milnor linking numbers \cite{Milnor54, Milnor57}. All of them, with an exception of Woltjer's helicity, are defined under restrictive assumptions either on the vector field $B$ or its domain. The reader should consult \cite{Khesin98, Khesin05, Moffatt95} for background material 
on helicity invariants and more specifically to Open Problem 7.18 posed by Arnold and Khesin in \cite[p. 176]{Khesin98} which asks to what extent these restrictive assumptions can be removed (see Remark \ref{rem:higher-helicities}).  
In physics, helicity invariants find a direct application in the phenomena of magnetic relaxation.
  An interested reader will find a thorough exposition of the subject in the work of Moffatt \cite{Moffatt85}. Specifically, Moffatt discusses why Borromean configurations of invariant tubes are relevant for the magnetic relaxation process. In short, if we minimize energy subject to keeping Woltjer's helicity constant, as Woltjer did in 1958 \cite{Woltjer58}, we obtain a ``force-free'' field (i.e. the {\em Beltrami field}). However, this field is not in fact realized under natural evolution, because its helicity is not the only invariant. Therefore, a construction of higher helicities may contribute to a better understanding of nature of the energy minimizers. The author is not aware however if 
 Borromean configurations have been observed in dynamical systems occurring in nature such as the magnetic fields on the Sun.

Throughout the article we use the convenient language of differential forms. In Section \ref{sec:1} we state our main results Theorem \ref{th:main} and Theorem \ref{th:main2} together with the necessary background. Section \ref{sec:2} is devoted to proofs of main theorems which use an integral formula for $\bar{\mu}_{123}$-invariant of $3$-component Borromean links in $\R^3$, this formula is stated in Key Lemma of the paper. A self contained exposition of all necessary background for Key Lemma and its proof are presented in Section \ref{sec:4} and the appendix.
 
\bigskip

\no \emph{Acknowledgments:} I wish to thank Professor Fred Cohen for constant support and for teaching me about configuration spaces,  I am equally grateful to Professor Paul Melvin for conversations about $\bar{\mu}$-invariants.

\section{Statement of results}\label{sec:2}

 Denote a parametrized $n$-component link in $\R^3$ (or $S^3$) by
$L=\{L_1,L_2,\ldots, L_n\}$, (where $L_i:S^1\mapsto \R^3 \text{(or $S^3$)}$ such that $L_i(S^1)\cap L_j(S^1)=\text{\O}$, $i\neq j$).  Recall that {\em link homotopy} is a deformation of a link which allows each component to pass through itself but not through a different component. The Milnor linking numbers also known as $\bar{\mu}$-invariants are invariants of $n$-component links up to link homotopy, we refer the reader to cf. \cite{Milnor57} for their definition.  Here, we will work entirely in the realm of $2$ or $3$-component links. For a $2$-component link $L=\{L_1,L_2\}$ there is just one $\bar{\mu}$-invariant, i.e. the linking number $\bar{\mu}_{12}(L_1,L_2)$ (or $\bar{\mu}_{12}$ when $L$ is known). In the language of intersection theory $\bar{\mu}_{12}$ is defined as the intersection number of one of the components of $L$ with a Seifert surface spanning the second component. Equivalently, we may define the linking number as the degree of 
a map from a 2-torus to the configuration space of two points in $\R^3$ (see Equation \eqref{eq:linking-gauss} and the discussion afterwards). It is well known that $\bar{\mu}_{12}$ is a complete invariant of $2$-component links up to link homotopy \cite{Milnor54}.  For $3$-component links the complete set of link homotopy invariants consists of the pairwise linking numbers $\bar{\mu}_{12}(L_1,L_2)$, $\bar{\mu}_{12}(L_2,L_3)$, $\bar{\mu}_{12}(L_1,L_3)$, and the triple linking number $\bar{\mu}_{123}\equiv\bar{\mu}_{123}(L_1,L_2,L_3)$ as an element  $\mathbb{Z}_{\text{gcd}(\bar{\mu}_{12}(L_1,L_2), \bar{\mu}_{12}(L_2,L_3), \bar{\mu}_{12}(L_1,L_3))}$ expressed in terms of the lower central series of the link group $G=\pi_1(S^3-L)$, cf. \cite{Milnor54}. In \cite{Mellor03} Mellor and Melvin found a geometric reformulation of Milnor's definition as follows: Choose Seifert surfaces $F_1$, $F_2$ and $F_3$ for the
components of $L=\{L_1, L_2, L_3\}$ and move these into general position. Starting at any point on $L_1$,record its intersection with the Seifert surfaces for $L_2$ and $L_3$ by a word $w_1$ in $2$ and $3$. For
example a $2$ or $2^{-1}$ in $w_1$ indicates a positive or negative intersection point of $L_1$ with $F_2$. Set $m_1 := m_{23}(w_1)$ to be a signed number of occurrences of $2$ and $3$ in the word $w_1$, for instance $\cdots 2\cdots 3\cdots$ or $\cdots 2^{-1}\cdots 3^{-1}\cdots$ contribute $+1$ to $m_1$, while $\cdots 2^{-1}\cdots 3\cdots$ or $\cdots 2\cdots 3^{-1}\cdots$ contribute $-1$ to $m_1$. Similarly, we define 
$m_2:= m_{31}(w_2)$ and $m_3:= m_{12}(w_3)$. We also let $t$ be the signed count
of the number of triple points of intersection of the three Seifert surfaces. Then
the {\em triple linking number} $\bar{\mu}_{123}(L)$ equals \cite{Mellor03}
\[
\bar{\mu}_{123}(L_1,L_2,L_3) = (m_1 + m_2 + m_3 - t) \mod \gcd(\bar{\mu}_{12}(L_1,L_2), \bar{\mu}_{13}(L_1,L_3), \bar{\mu}_{23}(L_2,L_3)).
\]
Note that if $L$ is Borromean, i.e. $\bar{\mu}_{12}(L_i,L_j)=0$, $i\neq j$, the  triple linking number is an integer valued invariant.

 So far the intersection theory approach to $\bar{\mu}$-invariants and their Massey product interpretation \cite{Porter80} was the main source of formulas for higher helicities  cf. \cite{Berger90, Akhmetiev05, Hornig04, Laurence-Stredulinsky00b} and \cite{Khesin98} for an overview. Here, we extend the methodology developed in \cite{Kom-unlinked09} based on the interpretation of $\bar{\mu}$-invariants as homotopy invariants 
 of associated link maps (see \cite{Koschorke04, Koschorke91, Koschorke97} and recently in \cite{Kom-Milnor08}). 

\begin{wrapfigure}{r}{0.46\textwidth}
	\begin{center}
   \includegraphics[width=.35\textwidth]{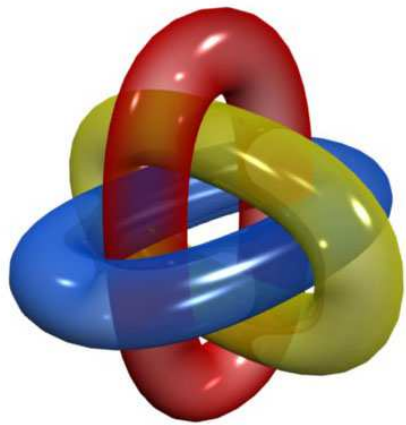}
  \end{center}
\caption{Flux tubes $\{\mathcal{T}_1,\mathcal{T}_2,\mathcal{T}_3\}$ modeled on the Borromean rings.}
	\label{fig:borro}
\end{wrapfigure} 

\no Let us denote by $(B,\mathcal{T}_\ast)$, a smooth vector field defined on the domain $\mathcal{T}_\ast$ which here we consider to  be either a closed manifold or a manifold with boundary, in the former case we additionally assume that $B$ is tangent to $\partial \mathcal{T}_\ast$. We will generally consider finitely many $B:=\{(B_i,\mathcal{T}_i)\}$, $i=1,\ldots, n$ where $\mathcal{T}_i$ are compact unless stated otherwise and let
$\mathcal{T}:=\prod^r_{i=1}\mathcal{T}_i$.

 Recall \cite{Vogel03, Khesin98} that a {\em system of short paths} on $\mathcal{T}_\ast$ is a collection of curves $\mathcal{S}=\{\sigma(x,y)\}$ indexed by pairs of points $(x,y) \in \mathcal{T}_\ast\times \mathcal{T}_\ast$ such that for any pair $(x,y)$ there is a connecting curve $\sigma(x,y):I\mapsto \mathcal{T}_\ast$, $\sigma(0)=x$ and $\sigma(1)=y$, and the lengths of curves in $\mathcal{S}$ are bounded by a common constant.  Given $T>0$ we introduce the following notation for orbits(left) and the {\em closed up orbits}(right) of a given $B_i$ after time $T$:
\begin{equation}\label{eq:orbits}
\begin{split}
 \mathscr{O}^{B_i}_T(x) & =  \{\Phi^i(x,t)\ |\ 0\leq t\leq T\},\\
 \bar{\mathscr{O}}^{B_i}_T(x)  & :=  \mathscr{O}^{B_i}_T(x)\cup \sigma(x,\Phi^i(x,T)),\quad \sigma(x,\Phi^i(x,T))\in \mathcal{S}.
 \end{split}
\end{equation}
\no In order to better motivate the definition of $\mathsf{H}_{123}(B;\mathcal{T})$ we first review the classical Woltjer's helicity which is defined for a pair of volume preserving vector fields  $B=\{(B_1,\mathcal{T}_1), (B_2,\mathcal{T}_2)\}$, $n=2$. Here we denote Woltjer's helicity by $\mathsf{H}_{12}(B)\equiv\mathsf{H}_{12}(B;\mathcal{T})$, \cite{Woltjer58, Arnold86}. The reader should consult \cite{Woltjer58} for the original definition of $\mathsf{H}_{12}(B)$.
Arnold's Helicity Theorem \cite{Arnold86, Vogel03} implies that Woltjer's helicity is given by the following integral
\begin{equation}\label{eq:H12} \mathsf{H}_{12}(B;\mathcal{T})=\int_{\mathcal{T}_1\times\mathcal{T}_2}\Bigl(\lim_{T\to\infty} \frac{1}{T^2}\bar{\mu}_{12}\bigl(\bar{\mathscr{O}}^{B_1}_T(x),\bar{\mathscr{O}}^{B_2}_T(y)\bigr) \Bigr)\,\vol_1(x)\wedge\vol_2(y),
\end{equation}
where $\vol_i$ denote volume forms on each $\mathcal{T}_i$ factor of $\mathcal{T}=\mathcal{T}_1\times\mathcal{T}_2$, and the function defined by the time 
average under the integral is referred to as {\em as the asymptotic linking number function}. The quantity on the right hand side is known as the \emph{average asymptotic linking number} \cite{Arnold86, Khesin05} or {\em asymptotic linking number} for short. It is currently unknown \cite{Khesin98} if  
 $\mathsf{H}_{12}(B)$ can be sensibly defined for vector fields not preserving 
 the volume element, but we may certainly assume the formula in \eqref{eq:H12} as a general definition of $H_{12}(B)$. In a similar spirit we 
define the third order helicity as an average asymptotic Milnor $\bar{\mu}_{123}$-invariant of orbits for triples $\{(B_i,\mathcal{T}_i)\}_{i=1,2,3}$, $n=3$. 

\begin{definition}\label{def:H123-definition}
Let $B:=\{(B_i,\mathcal{T}_i)\}_{i=1,2,3}$, be a triple of smooth vector fields defined above, then the third order helicity $\mathsf{H}_{123}(B;\mathcal{T})$ of $B$  is given by 
\begin{equation}\label{eq:H123-average}
 \mathsf{H}_{123}(B;\mathcal{T}):=\int_{\mathcal{T}} \Bigl(\lim_{T\to\infty} \frac{1}{T^3}\bar{\mu}_{123}\bigl(\bar{\mathscr{O}}^{B_1}_T(x),\bar{\mathscr{O}}^{B_2}_T(y),\bar{\mathscr{O}}^{B_3}_T(z)\bigr) \Bigr)\,\vol_1(x)\wedge\vol_2(y)\wedge\vol_3(z),
\end{equation}
whenever the limit under the integral 
\[
 \bar{m}_B:\,(x,y,z)\mapsto \lim_{T\to\infty} \frac{1}{T^3}\bar{\mu}_{123}\bigl(\bar{\mathscr{O}}^{B_1}_T(x),\bar{\mathscr{O}}^{B_2}_T(y),\bar{\mathscr{O}}^{B_3}_T(y)\bigr),
\]
exists almost everywhere and defines an integrable function $\bar{m}_B:\mathcal{T}\mapsto \R$ on $\mathcal{T}=\prod^3_{i=1} \mathcal{T}_i$ independent of the short paths system $\mathcal{S}$ chosen. Here $\vol_i$ denotes a volume form on the  $\mathcal{T}_i$ factor of $\mathcal{T}$. Subsequently, we refer to the function $\bar{m}_B$
 as the asymptotic $\bar{\mu}_{123}$-invariant function.
\end{definition}

\begin{figure}[htbp]
	\centering
	\begin{picture}(320, 210)
   \put(0,0){\includegraphics[width=.7\textwidth]{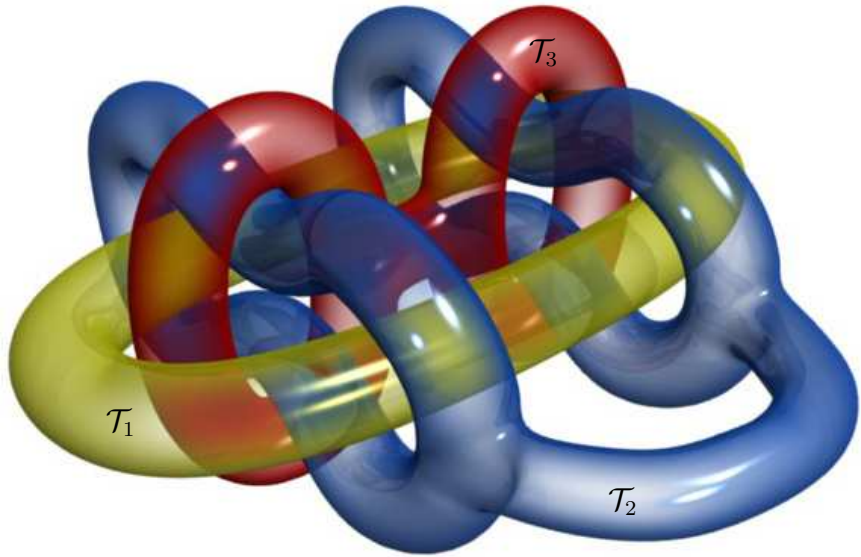}}
   \put(40,50){$\mathcal{T}_1$}
   \put(230,20){$\mathcal{T}_2$}
   \put(200,190){$\mathcal{T}_3$}
  \end{picture} 
\caption{Example of unlinked handlebodies $\mathcal{T}_1\cup\mathcal{T}_2\cup\mathcal{T}_3$.}
	\label{fig:invariant-sets}
\end{figure} 

In principle,  the above definition extends to higher Milnor linking numbers or generally to other link isotopy invariants \cite{Verjovsky94,Gambaudo-Ghys97,Baader-Marche07}. 

The main question one needs to address in the above definition is existence of the integral.  In \cite{Laurence-Stredulinsky00b} Laurence and Stredulinsky show existence of the third order helicity for {\em Borromean flux tubes}, i.e. domains $\mathcal{T}_i$ which are disjoint solid tori with cores forming a $3$-component Borromean link such as well known {\em Borromean rings} pictured on Figure \ref{fig:borro}. This type of domains are often referred to as domains \emph{modeled on a link}.  
 The main theorem of the current paper Theorem \ref{th:main} shows that $\mathsf{H}_{123}(B;\mathcal{T})$ is defined on unlinked invariant domains of $\R^3$ such as handlebodies pictured on Figure \ref{fig:invariant-sets}, and more importantly introduces a new formula for $\mathsf{H}_{123}(B;\mathcal{T})$. Before we state the main result we need to review several definitions. Recall 
\[
\Cnf_3(\R^3)=\{(x_1,x_2,x_3)\in \R^3\times\R^3\times\R^3;\ x_i\neq x_j, i\neq j\}.
\]
 Let $x=(u,v,r)\in \R^3\setminus \{0\}$, we define a closed differential 2-form
\begin{equation}\label{eq:omega}
 \omega(x):=\frac{1}{4\pi}\frac{u\,dv\wedge d r+v\,dr\wedge du+r\,du\wedge dv}{|x|^3},
\end{equation}
which restricts to the area form on the unit sphere in $S^2\subset \R^3$, normalized so that $\int_{S^2} \omega(x)=1$.
Define the \emph{Green form}
$\omega_{i,j}$ by
\begin{equation}\label{eq:w-ij}
 \omega_{i,j}:=\omega(x_i-x_j),\quad i>j.
\end{equation}
 In the vector notation 
\[ \omega_{i,j}(x_1,x_2,x_3)(X,Y)=\frac{\langle x_i-x_j, X, Y\rangle}{|x_i-x_j|^3},\qquad X,Y\in T(\R^3)^3,
\]
where $\langle\,\cdot\, ,\,\cdot\, ,\cdot\,\,\rangle$ denotes the triple product in $\R^3$.  It is well known \cite{Cohen-Lada-May76} that Green forms represent generators of the cohomology $H^\ast(\Cnf_3(\R^3))$ of the configuration space $\Cnf_3(\R^3)$ of three points in $\R^3$. (In Section \ref{sec:4}, we provide necessary background on the configuration space $\Cnf_3(\R^3)$.)

\begin{theorem}\label{th:main}
Suppose $\{\mathcal{T}_i\}_{i=1,2,3}$ are pairwise disjoint, compact, solid handlebodies $\mathcal{T}_i\subset \R^3$ with smooth boundary such  that every pair of 1-cycles  in $H_1(\mathcal{T}_i)$ and $H_1(\mathcal{T}_j)$, $i\neq j$ has linking number zero. Let $\{(B_i,\mathcal{T}_i)\}_{i=1,2,3}$ be a triple of volume preserving vector fields as defined above. Consider the integral 
\begin{equation}\label{eq:H123-integral-1}
\begin{split}
 \mathcal{J} & :=\int_{\mathcal{T}} \bigl(\omega_{1,2}\wedge d^{-1}\omega_{2,3} +\omega_{2,3}\wedge d^{-1}\omega_{3,1}+\omega_{3,1}\wedge d^{-1}\omega_{1,2}-\phi_{1,2,3}\bigr)\wedge \bigwedge^3_{i=1}\iota_{B_i}\vol_i,\\
 & \phi_{123}=d^{-1}(\omega_{1,2}\wedge\omega_{2,3}+\omega_{2,3}\wedge\omega_{3,1}+\omega_{3,1}\wedge\omega_{1,2}),\quad \iota_{B_i}\vol_i:=\vol_i(B_i,\,\cdot\,,\,\cdot\,)\ .
 \end{split}
\end{equation}
 Here $\vol_i$'s are volume forms of each $\R^3$ factor of $(\R^3)^3$.
Then,
 \begin{itemize}
 \item[$(i)$] $\mathsf{H}_{123}(B;\mathcal{T})$ exists and equals to $\mathcal{J}$. 
 
 \item[$(ii)$] $\mathcal{J}$ is  invariant of vector fields $B_i$ under the action of $\text{\rm SDiff}_0(\R^3)$, and an invariant of 2-forms $\iota_{B_i}\nu$ under the action of $\text{\rm Diff}_0(\R^3)$.
 
\end{itemize}
\end{theorem}

\begin{remark}
The smooth boundary assumption in the above theorem is in general not necessary cf. \cite{Kom-unlinked09}.
\end{remark}

\no In the following we will often abbreviate the sums under the integrals in Equation \eqref{eq:H123-integral-1} writing $\sum^3_{i=1} \omega_{i,i+1}\wedge \omega_{i+1,i+2}$ or $\sum^3_{i=1} \omega_{i,i+1}\wedge d^{-1} \omega_{i+1,i+2}$ understanding that the indices are taken modulo $3$.
\begin{remark}\label{rem:potentials}
Since the domain $\mathcal{T}$ of integration in \eqref{eq:H123-integral-1} is assumed to be  a product of disjoint compact solid handlebodies $\mathcal{T}_i\subset \R^3$ (as on Figure \ref{fig:invariant-sets}) such that every pair of 1-cycles $\mathscr{C}$,$\mathscr{C}'$ in 
$H_1(\mathcal{T}_i)$ and $H_1(\mathcal{T}_j)$, $i\neq j$ has a linking number zero
i.e.
\[
 \int_{\mathscr{C}\times\mathscr{C}'} \omega_{i,j}=0\ .
\]
Poincare duality implies that the 2-forms $\omega_{i,j}$ are exact on $\mathcal{T}=\prod_i\mathcal{T}_i$.
 Note that $\mathcal{J}$ is independent of a choice of 1-forms $\eta_{ij}\in d^{-1}\omega_{i,j}$ which we refer to as {\em potentials}  of $\omega_{i,j}$. Indeed, denoting 
\begin{equation}\label{eq:omega123-integral}
\begin{split}
 \omega_{1,2,3} & = \omega_{1,2}\wedge \eta_{2,3} +\omega_{2,3}\wedge \eta_{3,1}+\omega_{3,1}\wedge \eta_{1,2}-\phi_{1,2,3},\\
 d\eta_{i,j} & = \omega_{i,j},\quad \text{on}\quad \mathcal{T}_i\times \mathcal{T}_j,\quad i\neq j,\qquad \iota_B\vol:=\bigwedge^3_{i=1}\iota_{B_i}\vol_i\ .
\end{split}
\end{equation}
Observe that for any two potentials  $\eta_{i,j}$, $\eta'_{i,j}$ of $\omega_{i,j}$, the difference $\beta_{i,j}=\eta_{i,j}-\eta'_{i,j}$ is a closed $1$-form on $\mathcal{T}$ and therefore
\[
\omega_{1,2,3}-\omega'_{1,2,3} =\sum^3_{i=1}\omega_{i,i+1}\wedge (\eta_{i+1,i+2}-\eta'_{i+1,i+2})=d\Bigl(\sum^3_{i=1}\eta_{i,i+1}\wedge (\eta_{i+1,i+2}-\eta'_{i+1,i+2})\Bigr)=: d\beta\ .
\]
By Stokes Theorem
\[
 \mathsf{H}_{123}(B;\mathcal{T})-\mathsf{H}_{123}'(B;\mathcal{T}) = \int_{\mathcal{T}} (\omega_{1,2,3}-\omega'_{1,2,3})\wedge\iota_B\vol= \int_{\mathcal{T}} d\bigl(\beta\wedge \iota_B\vol\bigr)
  = \int_{\partial\mathcal{T}} \beta\wedge \iota_B\vol= 0,
\]
where in the second identity we applied $d(\iota_{B_i}\vol_i)=0$ as $B_i$'s are divergence free, and 
$\iota_{B_i}\vol_i\bigl|_{\partial\mathcal{T}_i}=0$ as each vector field $B_i$ is tangent to the boundary $\partial \mathcal{T}_i$.
\end{remark}

\no A crucial ingredient in the proof of Theorem \ref{th:main} (as may be expected from Definition \ref{def:H123-definition}) is the following

\begin{klemma}\label{lem:key-lemma}
{\em Given a parametrized $3$-component Borromean link $L=\{L_1$, $L_2$, $L_3\}$ in $\R^3$ denote by $F_L$ the associated product map
\[
 L_1\times L_2\times L_3:S^1\times S^1\times S^1\longrightarrow \Cnf_3(\R^3).
\]
Then, 
\begin{equation}\label{eq:I-mu-lem}
\begin{split}
  \bar{\mu}_{123}(L) & =\pm \int_{(S^1)^3} \bigl(F^\ast_L\omega_{1,2}\wedge \eta_{2,3} +F^\ast_L\omega_{2,3}\wedge \eta_{3,1}+F^\ast_L\omega_{3,1}\wedge \eta_{1,2}-F^\ast_L\phi_{1,2,3}\bigr),
 \end{split}
 \end{equation}
where $\eta_{i,j}$ satisfy $d\eta_{i,j}=F^\ast_L\omega_{i,j}$.}
\end{klemma}

\no A proof of Key Lemma will occupy Section \ref{sec:4} and Appendix.

\begin{remark} 
 Applying standard identities for Chen iterated integrals one may show that \eqref{eq:I-mu-lem} is equivalent to Chen's iterated integral proposed by Kohno in \cite[p. 155]{Kohno02}. 
\end{remark}

\begin{remark}[Explanation of terminology]
 Note that the flux formula for $\mathsf{H}_{123}(B;\mathcal{T})$, derived in \cite{Laurence-Stredulinsky00b, Kom-unlinked09} 
 shows it is $3^\text{rd}$-order in fluxes, in comparison, Woltjer's helicity $\mathsf{H}_{12}(B;\mathcal{T})$ is $2^\text{nd}$-order in fluxes \cite{Cantarella00}.
\end{remark}

\begin{remark}\label{rem:higher-helicities}
Woltjer's helicity $\mathsf{H}_{12}(B)$  is a well defined invariant of volume preserving vector fields $B_i$ defined on possibly ``overlapping'' domains. For instance 
we obtain the {\em self-helicity} $\mathsf{H}_{12}(\bar{B})$ when $(B_1,\mathcal{T}_1)=(B_2,\mathcal{T}_2)=(\bar{B},\bar{\mathcal{T}})$, and $\mathcal{\bar{T}}$ possibly a closed $3$-manifold such as $S^3$ or a homology sphere (otherwise the linking number is not well defined).  Open Problem 7.18 posed by Arnold and Khesin in \cite[p. 176]{Khesin98} asks if higher order invariants, such as $\mathsf{H}_{123}(B;\mathcal{T})$, can be defined possibly on overlapping domains. One may refer to these hypothetic invariants as {\em higher order asymptotic self-linking numbers}.  A recent work of Badder \emph{et al.} \cite{Baader07, Baader-Marche07} shows that for ergodic vector fields $B$ all asymptotic Vassiliev invariants are proportional to Wojtier's helicity, this may explain why attempts to define higher order asymptotic self-linking numbers have failed so far.  
\end{remark}

In magnetohydrodynamics cf. \cite{Priest84} magnetic fields evolve under the motion of supporting plasma (i.e. along a path in $\text{\rm SDiff}_0(\R^3)$). During the evolution, a magnetic field often dissipates its $L^2$-energy and among questions of interest
is whether the energy can be reduced to zero in the process cf. 
\cite{Khesin98}. Lower bounds for the $L^2$-energy of $B$ in terms of quantities, such as $\mathsf{H}_{123}(B;\mathcal{T})$ or $\mathsf{H}_{12}(B;\mathcal{T})$, invariant under the action of $\text{\rm SDiff}_0(\R^3)$ provide a way to decide this question for a given magnetic field $B$.
Woltjer's helicity provides such an energy bound, which is extensively used in e.g. magnetohydrodynamics, 
consult \cite{Khesin98} for further discussion. The next theorem applies formula \eqref{eq:H123-integral-1} to derive  a lower bound for the $L^2$-energy of $B$ in terms of  $\mathsf{H}_{123}(B;\mathcal{T})$. 

 Recall that the {\em Neumann Laplacian} $\Delta_N$ is the differential operator associated with the following boundary value problem \cite{Schwarz95}
 \begin{equation}\label{eq:neumann-problem}
 \begin{cases}
  \Delta \phi=\omega,\qquad \text{in}\quad \mathcal{T},\\
  \mathbf{n}\,\phi=\mathbf{n}\,d\phi=0,\qquad \text{on}\quad \partial \mathcal{T},
 \end{cases}
\end{equation}
where $\phi$, $\omega$ are differential forms of a fixed degree and $\mathbf{n}$ extracts the normal component of the form, $\Delta=d\delta+\delta d$ and $\delta=\pm\ast d\ast$. In other words $\Delta_N$ is the standard Laplace-Beltrami operator acting on the space of differential forms with boundary conditions specified in the above problem.

\begin{theorem}\label{th:main2}
 Let $\bar{B}$ be a volume preserving vector field in $\R^3$, given a triple of compact pairwise disjoint unlinked handlebodies $\mathcal{T}_i$ in $\R^3$ suppose $\bar{B}$ is tangent to $\partial\mathcal{T}_i$ for each $i=1,2,3$. Consider $B:=(B_i,\mathcal{T}_i)$, $B_i=\bar{B}|_{\mathcal{T}_i}$ as defined above. Then, the $L^2$-energy $E_2(\bar{B})=\int_{\R^3} |\bar{B}|^2$ of $\bar{B}$ admits the following lower bound
\[
 E_2(\bar{B})\geq \Bigl(|\mathsf{H}_{123}(B;\mathcal{T})|\cdot \frac{C\,r^2_\mathcal{T}\sqrt{\lambda_{1,N}}}{\sum^3_{i=1}\bigl(\| \omega_{i,i-1}\wedge\omega_{i,i+1}\|_2+\|\omega_{i,i+1}\|_2\bigr)}\Bigr)^{\frac 32}.
\]
Here the $L^2$-norm $\|\,\cdot\,\|_2$ is taken over $\mathcal{T}=\prod^3_{i=1} \mathcal{T}_i$, $\lambda_{1,N}$ is the first eigenvalue of the Neumann Laplacian on $\mathcal{T}$, $C$ is a universal constant and $r_\mathcal{T}$ denotes a minimal distance between pairs of handlebodies $\mathcal{T}_i$ in $\R^3$.
\end{theorem}


\section{Proof of Theorem \ref{th:main} and Theorem \ref{th:main2}}\label{sec:3}

\subsection{Third order helicity}
\no The following theorem is fundamental for our considerations.
\begin{theorem}[$L^1$-Ergodic Theorem, \cite{Becker81}]\label{th:L^1-ergodic}
 Given a triple of volume preserving flows $B:=(B_i,\mathcal{T}_i)$ and a real valued $L^1$-function $F$, consider $\bar{F}:\mathcal{T}\mapsto \R$, $\mathcal{T}=\prod^3_{i=1} \mathcal{T}_i$ which is called a time average of $F$, and is defined as follows 
\[
\bar{F}(x,y,z)=\lim_{T\to \infty}\frac{1}{T^3}\int^T_0\int^T_0\int^T_0
 F(\Phi_1(x,s),\Phi_2(y,t),\Phi_3(z,u))\,d s\,d t\,d u,
\]
where $\Phi_i(\,\cdot\,,\,\cdot\,)$ denotes the flow of $B_i$. Then,
\begin{itemize}
 \item[$(a)$] $\bar{F}(x,y,z)$ exists almost everywhere,
 \item[$(b)$] $\|\bar{F}\|_{L^1(\mathcal{T})}\leq \|F\|_{L^1(\mathcal{T})}$,
 \item[$(c)$] $\bar{F}$ is invariant under the action by the flows of $B_i$,
 \item[$(d)$] if $\mathcal{T}$ is of finite volume then
\begin{gather}\label{eq:ergodic-integrals-equal}
 \int_{\mathcal{T}}\bar{F} =\int_{\mathcal{T}} F\ .
\end{gather}
\end{itemize}
\end{theorem}

\begin{proof}[Proof of Theorem \ref{th:main}]
We first show $(i)$, note that 
 the following identity proven in Appendix A of \cite{Kom-unlinked09} is valid for any 3-form $\beta$ on $\mathcal{T}$ and a triple of fields $B_i$ on $\mathcal{T}_i$:
 \begin{equation}\label{eq:iota-volume}
 \begin{split}
(\iota_{B_3}\iota_{B_2}\iota_{B_1}\beta)\wedge\vol_1\wedge\vol_2\wedge\vol_3 & =
\beta(B_1,B_2,B_3)\,\vol_1\wedge\vol_2\wedge\vol_3\\
 & =\beta\wedge\iota_{B_1}\vol_1\wedge\iota_{B_2}\vol_2\wedge\iota_{B_3}\vol_3\ .
 \end{split}
 \end{equation}
 As a result \begin{eqnarray*}
\notag  \mathcal{J} & = & \int_{\mathcal{T}} \omega_{1,2,3}\wedge \iota_{B_1}\vol_1\wedge \iota_{B_2}\vol_2\wedge \iota_{B_3}\vol_3 = \int_\mathcal{T} \bigl(\iota_{B_3}\iota_{B_2}\iota_{B_1}\omega_{1,2,3}\bigr)\,\vol_1\wedge\vol_2\wedge\vol_3\\
\label{eq:m-function}  & = & \int_\mathcal{T} m_B(x,y,z)\,\vol_1(x)\wedge \vol_2(y)\wedge \vol_3(z),
\end{eqnarray*}
where $m_B:=\omega_{1,2,3}(B_1,B_2,B_3)$. Thanks to the assumptions on $\mathcal{T}$, every triple of closed up orbits $\{\bar{\mathscr{O}}^{B_1}_T(x),\bar{\mathscr{O}}^{B_2}_T(y),\bar{\mathscr{O}}^{B_3}_T(y)\}$ is a $3$-component Borromean link. By Key Lemma
\begin{eqnarray}
\notag\bar{\mu}_{123}(\bar{\mathscr{O}}^{B_1}_T(x),\bar{\mathscr{O}}^{B_2}_T(y),\bar{\mathscr{O}}^{B_3}_T(y)) & = &
\int_{\bar{\mathscr{O}}^{B_1}_T(x)\times\bar{\mathscr{O}}^{B_2}_T(y)\times\bar{\mathscr{O}}^{B_3}_T(y)}
\omega_{1,2,3}\\
\label{eq:mu-O} & = & \int^T_0 \int^T_0 \int^T_0 \omega_{1,2,3}(B_1,B_2,B_3)\,ds\, dt\, du+(I)\\
\notag & = & \int^T_0\int^T_0\int^T_0 m_B(\Phi^1(x,s),\Phi^2(y,t),\Phi^3(z,u))\,ds\, dt\, du+(I),
\end{eqnarray}
where $(I)$ denotes integrals over short paths. 
\begin{remark}
 As a consequence of 
 homotopy invariance proven in Proposition \ref{th:integral-whitehead} and the fact that any piecewise smooth link 
 may be approximated by a smooth one, the integral \eqref{eq:I-mu} is well defined for piecewise smooth links such as $\{\bar{\mathscr{O}}^{B_1}_T(x),\bar{\mathscr{O}}^{B_2}_T(y),\bar{\mathscr{O}}^{B_3}_T(y)\}$. 
\end{remark}
\no Now, since short paths have bounded length we obtain
\begin{equation}\label{eq:I-to-zero}
\frac{1}{T^3} (I)\longrightarrow 0,\quad \text{as}\quad T\to\infty\ .
\end{equation}
Therefore, the expression in \eqref{eq:mu-O} for  $\bar{\mu}_{123}(\bar{\mathscr{O}}^{B_1}_T(x),\bar{\mathscr{O}}^{B_2}_T(y),\bar{\mathscr{O}}^{B_3}_T(y))$ combined with Theorem \ref{th:L^1-ergodic} $(a)$, $(d)$ and Definition \ref{def:H123-definition} yields
\[
  \mathcal{J}=\int_\mathcal{T} m_B=\int_\mathcal{T} \bar{m}_B=\mathsf{H}_{123}(B;\mathcal{T}),
\]
where $\bar{m}_B$ is the time average of $m_B$ which we called in Definition \ref{def:H123-definition} the asymptotic $\bar{\mu}_{123}$-function. Observe that, 
thanks to \eqref{eq:I-to-zero}, $\bar{m}_{B}$ is independent on the short path system chosen, thus we verified  existence of $\mathsf{H}_{123}(B;\mathcal{T})$ in the assumed setting, as well as the formula $\mathsf{H}_{123}(B;\mathcal{T})=\mathcal{J}$.

\no The proof of $(ii)$ is in the style of 
\cite{Berger90, Mayer03}, but adapted to our setting. 
For any given $g\in \text{\rm SDiff}(\R^3)$, by definition, there exists a path $t\longrightarrow g(t)\in \text{\rm SDiff}_0(\R^3)$,
such that
\[
 g(0)=\text{id}_{\R^3},\qquad g(1)=g\ .
\]
 Denote by $V$ the divergence free vector field on $\R^3$, given by $V(x)=\frac{d}{dt} g(t,x)|_{t=0}$, i.e. $g(t)$ is a flow of $V$. Let the push-forward fields $B_i$ be
\begin{equation}\label{eq:vol-pres-action}
 B^t_i:=g(t)_\ast B_i\ .
\end{equation}
It is well known \cite[p. 224]{Freedman91-2}
 that 2-forms: $\iota_{B^t_i}\vol$ are frozen in the flow of $V$ i.e.
\begin{equation}\label{eq:forms-frozen-in}
 \frac{d}{dt}\bigl(g(t)^\ast\iota_{B^t_i}\vol\bigr)|_{t=0} =(\partial_t+\mathcal{L}_V) \iota_{B^t_i}\vol_i = 0.
\end{equation} 
The tangent bundle $T(\R^3)^3$ has a natural product structure and we also have the path $\hat{g}(t)=(g(t),g(t),g(t))$ in $\text{\rm SDiff}_0(\R^3\times \R^3\times \R^3)$, which leads to the vector field $\hat{V}=(V,V,V)=\frac{d}{dt} \hat{g}(t)$. Equation \eqref{eq:forms-frozen-in} implies
\begin{equation}\label{eq:forms-frozen-in2}
(\partial_t+\mathcal{L}_{\hat{V}}) \iota_{B^t}\vol= 0,\qquad \text{where}\quad\iota_{B^t}\vol:=\iota_{B^t_1}\vol_1\wedge \iota_{B^t_2}\vol_2\wedge \iota_{B^t_3}\vol_3\ .
\end{equation} 

\begin{remark}\label{rem:diff}
 Notice that if we consider the action of $\text{\rm Diff}_0(\R^3)$ on forms $\iota_{B_i}\vol_i$ by pullbacks:
 $(g,\iota_{B_i}\vol_i)\longrightarrow (g^\ast)^{-1}(\iota_{B_i}\vol_i)$, and let $(\iota_B\vol)^t:=(g^\ast(t))^{-1}\iota_B\vol$
  we immediately obtain $(\partial_t+\mathcal{L}_{\hat{V}}) (\iota_B\vol)^t = \mathcal{L}_{\hat{V}} (\iota_B\vol)^t=0$.
\end{remark}

\no Let $\mathcal{T}(t):=\hat{g}(t)(\mathcal{T})\subset \Cnf_3(\R^3)$, we must show $\frac{d}{dt} 
 \mathsf{H}_{123}(g(t)_\ast B ,\mathcal{T}(t))=0$. 
Without loss of generality we set 
$t_0=0$, since for any $t=t_0$ we may apply a pullback by $\hat{g}(t_0)$. Then $\hat{g}(0)=\text{id}_{(\R^3)^3}$ and we have
\begin{equation}
\begin{split}
 \frac{d}{dt} 
 \Bigl(\mathsf{H}_{123}(g(t)_\ast B ,\mathcal{T}(t))\Bigr)\Bigl|_{t=0} & = 
 \frac{d}{dt} \int_{\mathcal{T}(t)} \omega_{1,2,3}\wedge \iota_{B^t}\vol= \int_{\mathcal{T}(0)} \frac{d}{dt} \hat{g}(t)^\ast\bigl(\omega_{1,2,3}\wedge \iota_{B^t}\vol\bigr)\\
 & =  \int_{\mathcal{T}(0)} \bigl(\mathcal{L}_{\partial_t+\hat{V}}(\omega_{1,2,3})\bigr)\wedge \iota_{B^t}\vol,
\end{split}
\end{equation}
where in the last identity we applied \eqref{eq:forms-frozen-in2} and the product rule for the Lie derivative.
Because $\omega_{1,2,3}$ is time independent, and $d\omega_{1,2,3}=0$, Cartan's magic formula yields
\begin{eqnarray*}
 \mathcal{L}_{\partial_t+\hat{V}}(\omega_{1,2,3}) & = & \mathcal{L}_{\hat{V}}(\omega_{1,2,3}) = d(\iota_{\hat{V}} (\omega_{1,2,3}))\ .
\end{eqnarray*}
Since $B^t_i$ are tangent to the boundary of $\mathcal{T}_i(t)$ the same argument as in Remark \ref{rem:potentials} shows that the right hand side of the previous equation vanishes. Thanks to Remark \ref{rem:diff}, we  obtain the second statement of $(ii)$ analogously.
\end{proof}

\subsection{Lower bound for the \texorpdfstring{$L^2$}{L2}-energy}

\begin{proof}[Proof of Theorem \ref{th:main2}]
From the Cauchy inequality (where $\ast$ is the Hodge star operator) and derivation in \cite[p. 22]{Kom-unlinked09},
 we estimate
\[
 |\mathsf{H}_{123}(B;\mathcal{T})|=|\langle\ast\omega_{1,2,3}, \iota_B\vol\rangle_2|  \leq \|\omega_{1,2,3}\|_2\|\iota_B\vol\|_2\leq \|\omega_{1,2,3}\|_2\, E_2(B)^{\frac{3}{2}},
\]
 The norm $\|\omega_{1,2,3}\|_2$ can be bounded using geometry of $\R^3$ (rather than 
the round metric of $S^3$ as in \cite{Kom-unlinked09}) as follows. Let $d^{-1}:=\delta\Delta^{-1}_N$, where $\Delta_N$ is 
the Neumann Laplacian on differential 1-forms on $\mathcal{T}\subset (\R^3)^3$, then the $3$-form $\phi_{123}$ in \eqref{eq:H123-integral-1} is given by  $\phi_{123}:=d^{-1}\bigl(\sum^3_{i=1} \omega_{i,i+1}\wedge \omega_{i+1,i+2}\bigr)$. We obtain
\[
\begin{split}
 \|\omega_{1,2,3}\|_2 &=\|\omega_{1,2}\wedge d^{-1}\omega_{2,3} +\omega_{2,3}\wedge d^{-1}\omega_{3,1}+\omega_{3,1}\wedge d^{-1}\omega_{1,2}-\phi_{123}\|_2\\
 & \leq \|d^{-1}\|\|\omega_{r_{\mathcal{T}}}\|_\infty\Bigl(\sum^3_{i=1}\bigl(\| \omega_{i,i-1}\wedge\omega_{i,i+1}\|_2+\|\omega_{i,i+1}\|_2\bigr)\Bigr),
\end{split}
\]
where $\omega_r=\omega|_{\R^3-B(r)}$ denotes the restriction of $\omega$ defined in \eqref{eq:omega} to the complement of a
radius $r$ ball $B(r)\subset \R^3$, and $r_\mathcal{T}$ is a lower bound for the minimum distance between pairs of handlebodies $\mathcal{T}_i$ in $\R^3$. Clearly, $\|\omega_r\|_\infty$ grows like $\frac{1}{r^2}$ as $r\to 0$ thus $\|\omega_{r_\mathcal{T}}\|_\infty\geq C^{-1}\,r^{-2}_{\mathcal{T}}$ for some universal constant $C$.
 Since $\|d^{-1}\|\leq \frac{1}{\sqrt{\lambda_{1,N}}}$ (c.f \cite{Schwarz95}) where $\lambda_{1,N}$ is the first eigenvalue of the Neumann Laplacian on $\mathcal{T}$, we obtain the estimate as claimed.
\end{proof} 

We expect that the presented method will lead to a hierarchy of helicities defined 
on invariant $n$-component unlinked domains together with associated 
energy bounds \cite{Kom-hierarchy}.


\section{Integral formula for Milnor \texorpdfstring{$\bar{\mu}_{123}$}{mu123}-invariant.} \label{sec:4}

\subsection{Background on \texorpdfstring{$\Cnf_3(\R^3)$}{Conf3(R3)}.}

 Following \cite{Fadell-Husseini01} we set $e$ to be the unit vector $(1,0,0)$ in $\R^3$ and define
\[
 q_1=0,\quad q_2=4\,e,\quad q_3=8\,e,\quad Q_i=\{q_1,\ldots,q_i\},\quad Q_0=\text{\O}\ .
\]
\no The following spherical cycles on $\Cnf_3(\R^3)$ are of fundamental importance
\begin{equation}\label{eq:Aij-def}
\begin{split}
 A_{i,j} & : S^2\longrightarrow \Cnf_3(\R^3),\qquad 1\leq j < i \leq 3,\\
     A_{2,1} & : \xi \longrightarrow (q_1,\xi,q_3),\quad A_{3,2}  : \xi \longrightarrow (q_1,q_2,q_2+\xi),\quad 
     A_{3,1}  : \xi \longrightarrow (q_1,q_2,\xi).
\end{split}
\end{equation}
We denote their respective homotopy classes in $\pi_2(\Cnf(\R^3))$ by $\alpha_{i,j}$. Consider projections 
\begin{equation}\label{eq:Pi_i}
\begin{split}
\Pi_i: & \Cnf_3(\R^3)\longrightarrow \Cnf_2(\R^3),\\
\Pi_i & (x_1, x_2, x_3)  = (\ldots,\widehat{x}_i,\ldots),\quad i=1,2,3,
\end{split}
\end{equation}
defined by skipping the $i$-th coordinate factor. Because $\Cnf_2(\R^3)$ is diffeomorphic to 
$\R^3\times \bigl(\R^3-\{0\}\bigr)$, via $(x_1,x_2)\mapsto (x_1, x_2-x_1)$, it has a homotopy type of $S^2$. Directly from the 
definition it follows that  $\Pi_k\circ A_{i,j}$ are degree one maps when $i,j\neq k$ or 
null homotopic whenever $i=k$ or $j=k$. Results of \cite{Fadell-Husseini01,Cohen-Lada-May76} tell us that every $\Pi_i$ is a fibration which admits a section. In particular, choosing $i=3$ we obtain the fibration diagram
\begin{equation}\label{eq:fibration-pi}
  \CD
  \Cnf_3(\R^3) @<<<  \Cnf_1(\R^3-Q_2)=\R^3-Q_2\cong S_{3,1}\vee S_{3,2}\\
  @VV\Pi_3 V \\
  \Cnf_2(\R^3)\cong S_{2,1}
    \endCD
\end{equation}
\no where by $S_{i,j}$ we denote the images: $A_{i,j}(S^2)\subset \Cnf_3(\R^3)$. Obviously, we may choose to fiber over each $S_{i,j}$ separately. As an immediate consequence, 
we obtain \cite[p. 189]{Whitehead78}
\begin{equation}\label{eq:pi_k}
 \pi_k(\Cnf_3(\R^3))\cong \pi_k(S_{2,1})\oplus\pi_k(S_{3,1}\vee S_{3,2}).
\end{equation}
In particular for $k=2$, we conclude that $\alpha_{i,j}$ generate $\pi_2(\Cnf_3(\R^3))\cong \Z\oplus \Z\oplus \Z$. 
 Next, we describe a structure of the deRham cohomology ring of the configuration space $\Cnf_3(\R^3)$, \cite{Cohen-Lada-May76}. 
Every $\omega_{i,j}$ represents an integral cohomology class $\psi_{i,j}:=[\omega_{i,j}]$ and is dual to the cycle $[A_{i,j}]$ defined in \eqref{eq:Aij-def}. The cohomology ring $H^\ast_{\text{dR}}(\Cnf_n(\R^3))$ is generated \cite{Cohen-Lada-May76} by $\psi_{i,j}$, $1\leq j < i\leq 3$ with relations
\begin{equation}\label{eq:conf-cohomology-rel}
\begin{split}
 \psi^2_{i,j} & =0,\qquad \psi_{i,j}=-\psi_{j,i},\\
 \psi_{3,1}\psi_{3,2} & =\psi_{2,1}(\psi_{3,2}-\psi_{3,1}),
\end{split}
\end{equation}
see \cite{Cohen-Lada-May76} and \cite[p. 101]{Fadell-Husseini01}.
The last relation on representatives $\omega_{i,j}$ reads
\begin{equation}\label{eq:wij-relation}
\omega_{2,1}\wedge\omega_{3,2}-\omega_{3,2}\wedge\omega_{3,1}-\omega_{3,1}\wedge\omega_{2,1}=\sum^3_{i=1} \omega_{i,i+1}\wedge\omega_{i+1,i+2}=d\phi_{1,2,3},
\end{equation}
for some smooth 3-form $\phi_{1,2,3}$.


\subsection{Whitehead products in the configuration space \texorpdfstring{$\Cnf_3(\R^3)$}{Conf3(R3)}}

\no Our goal in a later section is to understand $\pi_3(\Cnf_3(\R^3))$ in the context of so called link maps. Thanks to the decomposition in \eqref{eq:pi_k} among relevant generators of this group are the Whitehead products  of $\alpha_{i,j}$'s  \cite{Whitehead78}. We aim to obtain suitable integrals for these Whitehead products.

  Let $D^{p}$ denote a $p$ dimensional disk in $\R^{p+1}$, given two continuous maps 
 \[
  f_k:(D^{p_k},\partial D^{p_k})\longrightarrow (X,x_0),\quad k=1,2,
\]  
into a pointed topological space $(X,x_0)$, the Whitehead product of $f_1$ and $f_2$ is given by \cite{Whitehead78}
\begin{equation}\label{eq:whitehead-product}
\begin{split}
 [f_1,f_2]:\partial (D^{p_1}\times  D^{p_2}) & \cong S^{p_1+p_2-1}\longrightarrow (X,x_0),\\
[f_1,f_2](x_1,x_2) & :=
\begin{cases}
 f_1(x_1),\qquad x_2\in  \partial D^{p_2},\\
 f_2(x_2),\qquad x_1\in  \partial D^{p_1},
\end{cases}
\end{split}
\end{equation}
recall $\partial (D^{p_1}\times  D^{p_2})=\partial D^{p_1}\times  D^{p_2}\cup D^{p_1}\times  \partial D^{p_2}$. The operation $[\,,\,]:\pi_{p_1}(X)\times \pi_{p_2}(X)\longrightarrow \pi_{p_1+p_2-1}(X)$ is well defined and turns the vector space $\pi_\ast(X)\otimes \R$ into a graded Lie algebra over $\R$ cf.\cite{Cohen-Lada-May76}. 
In the following proposition we extend calculations in \cite{Haefliger78} to define an integral
detecting certain Whitehead products in the configuration space $\Cnf_3(\R^3)$ (also compare with Section 3.3 in the preprint \cite{Sinha-Walter08}).

\begin{proposition}\label{th:integral-whitehead}
 For any $f:S^3\longrightarrow \Cnf_3(\R^3)$ let
 \begin{equation}\label{eq:I}
  \mathcal{I}(f):=\Bigl(\int_{S^3} \sum^3_{i=1}f^\ast\omega_{i,i+1}\wedge \eta_{i+1,i+2}\Bigr) -\int_{S^3} f^\ast\phi_{1,2,3},
 \end{equation}
 where forms $f^\ast\omega_{i,j}$ are exact and $d \eta_{i,j}=f^\ast\omega_{i,j}$.
 Then,
 \begin{itemize}
 \item[$(i)$] $\mathcal{I}$ is independent of the choice of potentials $\eta_{i,j}$.
 \item[$(ii)$]  $\mathcal{I}\in \text{\rm Hom}(\pi_3(\Cnf_3(\R^3)),\R)$ and satisfies
 \begin{equation}\label{eq:I-Aij}
  \begin{split}
  \mathcal{I}([\alpha_{1,2},\alpha_{2,3}]) & =\mathcal{I}([\alpha_{2,3},\alpha_{3,1}])=\mathcal{I}([\alpha_{3,1},\alpha_{1,2}])=1,\\
  \mathcal{I}([\alpha_{i,j},\alpha_{i,j}]) & =0,\quad 1\leq j < i\leq 3\ .
  \end{split}
 \end{equation}
 \end{itemize}
\end{proposition}
\begin{proof}
 $\mathcal{I}$ is independent of a choice of potentials $\eta_{i,j}$'s: indeed, let $\eta'_{ij}$ be different potentials then $d(\eta_{i,j}-\eta'_{i,j})=f^\ast\omega_{i,j}-f^\ast\omega_{i,j}=0$ and
 \[
  \begin{split}
  \mathcal{I}(f)-\mathcal{I}'(f) & =\sum^3_{i=1} \int_{S^3} f^\ast\omega_{i,i+1}\wedge (\eta_{i+1,i+2}-\eta'_{i+1,i+2}) =\sum_{i=1} \int_{S^3} d\eta_{i,i+1}\wedge (\eta_{i+1,i+2}-\eta'_{i+1,i+2})\\
  & =\sum^3_{i=1} \int_{S^3} \eta_{i,i+1}\wedge d(\eta_{i+1,i+2}-\eta'_{i+1,i+2})=0,
  \end{split}
 \] 
 where in the third identity we applied Stokes Theorem. To show that $\mathcal{I}$ is a well defined homomorphism we first show invariance under homotopies. Let $F:I\times S^3\to \Cnf_3(\R^3)$ be a
 homotopy between $F_0$ and $F_1$, and let $\tilde{\eta}_{i,j}=d^{-1}F^\ast\omega_{i,j}$ on $S^3\times I$. Combining Stokes Theorem, $(i)$, Equation \eqref{eq:wij-relation} and $d\eta_{i,j}=\omega_{i,j}$ we obtain
\begin{eqnarray}
   \label{eq:I-homotopy-invariance} \mathcal{I}(F_1)-\mathcal{I}(F_0) & = &
  \int_{S^3\times \{1\}} \bigl(\sum^3_{i=1} F^\ast_1\omega_{i,i+1}\wedge 	\eta_{i+1,i+2}-F^\ast_1\phi_{1,2,3}\bigr)\\
\notag  & &\quad -\int_{S^3\times \{0\}} \bigl(\sum^3_{i=1} F^\ast_0\omega_{i,i+1}\wedge 	\eta'_{i+1,i+2}-F^\ast_0\phi_{1,2,3}\bigr)\\
\notag  & = & \int_{S^3\times I} d\bigl(\sum^3_{i=1} F^\ast\omega_{i,i+1}\wedge \tilde{\eta}_{i+1,i+2}-F^\ast\phi_{1,2,3}\bigr)\\
\notag  & = & \int_{F(S^3\times I)} (\sum^3_{i=1} \omega_{i,i+1}\wedge d\tilde{\eta}_{i+1,i+2}-d\phi_{1,2,3}\bigr)=0\ .
\end{eqnarray}
 Additivity of $\mathcal{I}$ is a direct consequence of additivity for integrals and the definition of $+$ in $\pi_n(\,\cdot\,)$, thus $\mathcal{I}$ is a well defined element of $\text{\rm Hom}(\pi_3(\Cnf_3(\R^3)),\R)$.
 
  For Equations \eqref{eq:I-Aij},
consider $f=[f_1,f_2]$ defined in \eqref{eq:whitehead-product}. Let $p_1=p_2=2$, and $\pi_i:D_1\times D_2\mapsto D_i$
$i=1,2$ be projections onto each factor $D_i\cong D^2$, and $j:\partial(D_1\times D_2)\mapsto D_1\times D_2$ the inclusion. According to \eqref{eq:whitehead-product}, we have
\[
 f^\ast \omega_{i,j}=j^\ast(\pi^\ast_1 f^\ast_1\omega_{i,j}+\pi^\ast_2 f^\ast_2\omega_{i,j}).
\]
Because every 2-form $\pi^\ast_1 f^\ast_1\omega_{i,j}+\pi^\ast_2 f^\ast_2\omega_{i,j}$ is exact on $D_1\times D_2$, we may define a smooth potential $\eta_{i,j}$ such that
\begin{equation}\label{eq:d-eta}
 d\eta_{i,j}=\pi^\ast_1 f^\ast_1\omega_{i,j}+\pi^\ast_2 f^\ast_2\omega_{i,j},
\end{equation}
clearly $d j^\ast  \eta_{i,j}=f^\ast\omega_{i,j}$.
We calculate by applying Stokes Theorem 
and \eqref{eq:d-eta}
\[
\begin{split}
\mathcal{I}(f) & =\int_{\partial(D_1\times D_2)} \bigl(\sum^3_{i=1} f^\ast\omega_{i,i+1}\wedge \eta_{i+1,i+2}-f^\ast\phi_{1,2,3}\bigr)\\
& = \int_{D_1\times D_2}\bigl(\sum^3_{i=1}(\pi^\ast_1 f^\ast_1\omega_{i,i+1}+\pi^\ast_2 f^\ast_2\omega_{i,i+1})\wedge
(\pi^\ast_1 f^\ast_1\omega_{i+1,i+2}+\pi^\ast_2 f^\ast_2\omega_{i+1,i+2})\\
& \qquad-(\pi^\ast_1 f^\ast_1 d\phi_{1,2,3}+\pi^\ast_2 f^\ast_2 d\phi_{1,2,3})\bigr)\\
& =\int_{D_1\times D_2}\bigl(\sum(\pi^\ast_1 f^\ast_1(\omega_{i,i+1}\wedge\omega_{i+1,i+2})+\pi^\ast_2 f^\ast_2(\omega_{i,i+1}\wedge\omega_{i+1,i+2})\\
&\qquad -(\pi^\ast_1 f^\ast_1 d\phi_{1,2,3}+\pi^\ast_2 f^\ast_2 d\phi_{1,2,3})\bigr)\\
&\qquad  + \int_{D_1\times D_2}\bigl(\sum(\pi^\ast_1 f^\ast_1\omega_{i,i+1}\wedge \pi^\ast_2 f^\ast_2\omega_{i+1,i+2}+\pi^\ast_2 f^\ast_2\omega_{i,i+1}\wedge
\pi^\ast_1 f^\ast_1\omega_{i+1,i+2})\bigr)\ .
\end{split}
\]
\no The first integral in the above identity vanishes because of relations in \eqref{eq:wij-relation}. The second integral is equal to
\[
\begin{split}
 \mathcal{I}(f) & =\sum^3_{i=1}(\int_{D_1}f^\ast_1\omega_{i,i+1}\int_{D_2} f^\ast_2\omega_{i+1,i+2}+(-1)^4\int_{D_1}f^\ast_1\omega_{i+1,i+2}\int_{D_2}f^\ast_2\omega_{i,i+})\\
 & =\sum^3_{i=1} \bigl( \omega_{i,i+1}(f_1)\omega_{i+1,i+2}(f_2)+\omega_{i+1,i+2}(f_1)\omega_{i,i+1}(f_2)\bigr),
\end{split}
\]
where $\omega_{i,j}(f)=\int_{S^2} f^\ast\omega_{i,j}$. Identities in \eqref{eq:I-Aij} follow from
the definition of $A_{i,j}$ in \eqref{eq:Aij-def} and \eqref{eq:conf-cohomology-rel}.
\end{proof}

\begin{remark}
The above proposition certainly can be obtained from theories developed in \cite{Novikov88, Sullivan77, Hain84} and more recently 
in \cite{Sinha-Walter08}. Introduction of these theories would require a significant detour and is outside of the scope of this paper.
\end{remark}

\begin{remark}
Let $X$ be a smooth manifold, the argument in \cite{Haefliger78} ``runs'' as follows: let $\omega_1$ and $\omega_2$ be closed differential forms of degree $p_1$ and $p_2$, such that $\omega_1\wedge\omega_2=0$. Consider  spherical cycles $f_k:S^{p_k}\longrightarrow X$, $k=1,2$. It is shown that for any $[f]\in\pi_{p_1+p_2-1}(X)$
\[
 \mathcal{J}_{(\omega_1,\omega_2)}(f):=\int_{S^{p_1+p_2-1}} f^\ast\omega_1\wedge \eta_2=\int_{S^{p_1+p_2-1}} f^\ast\omega_2\wedge \eta_1,
\]
where $f^\ast\omega_k=d\eta_k$, $k=1,2$, defines an element of $\text{\rm Hom}(\pi_{p_1+p_2-1}(X);\R)$ satisfying
\begin{equation}\label{eq:whitehead-eval}
 \mathcal{J}_{(\omega_1,\omega_2)}([f_1,f_2])=\omega_1(f_1)\omega_2(f_2)+(-1)^{p_1 p_2}\omega_1(f_2)\omega_2(f_1),
\end{equation}
In particular given a degree one map, $f:S^p\to S^p$, Equation \eqref{eq:whitehead-eval} implies
 \begin{equation}\label{eq:2hopf}
  \mathcal{J}_{(\omega,\omega)}([f,f])=\begin{cases}
   0,\quad p=2k+1,\\
   2,\quad p=2k,
  \end{cases}\ .
\end{equation}
for a volume form $\omega$ on $X=S^p$, such that $\int_{S^p} \omega=1$.
Therefore, for even $p$, $[f,f]:S^{2p-1}\to S^p$ is twice the Hopf map and null for odd $p$.

\end{remark}


\subsection{Link maps.}

\no Denote a parametrized $n$-component link in $\R^3$ by
$L=\{L_1,L_2,\ldots, L_n\}$, (where $L_i:S^1\mapsto \R^3$, such that $L_i(S^1)\cap L_k(S^1)=\text{\O}$, $i\neq k$), $L$
defines a \emph{link map}, cf. \cite{Koschorke91, Koschorke97}:
\[
 L: \bigsqcup^n_{i=1} S^1\longrightarrow \R^3,\quad L\bigl|_{S^1_i}=L_i.
\]
We denote by  $LM(n)$ the set of link homotopy classes of $n$-component link maps. 
In  \cite{Koschorke91, Koschorke97} the author defines the $\kappa$-invariant
\begin{equation}\label{eq:link-map}
\begin{split}
    \kappa & : LM(n)\longrightarrow [(S^1)^n,\Cnf_n(\R^3)],\\
 \kappa(L)  = [ F_L & :(S^1)^n\longrightarrow \Cnf_n(\R^3) ],\quad  F_L=L_1\times\ldots\times L_n\ .
\end{split}
\end{equation}
\no $\kappa(L)$ is well defined because a link homotopy of $L$ in $\R^3$ yields a homotopy of the associated $F_L$. Note that the set of based homotopy classes is in bijective correspondence with the set of base point free homotopy classes  because  $\Cnf_n(\R^3)$ is simply connected cf. \cite{Whitehead78}. 
$\kappa$-invariants are closely tied to $\bar{\mu}$-invariants. It is has been proven by Koschorke in Corollary 6.2 of \cite[p. 314]{Koschorke97} that whenever $L$ is a Borromean $n$-component  
link, $\kappa(L)$ can be identified, up to a sign, with $(n-2)!$ integers which are all possible  $\bar{\mu}$-invariants of $L$.

 Let us review the basic case of the linking number $\bar{\mu}_{12}(L_1, L_2)$ in $\R^3$. Denote parametrizations of components by  $L_1=\{x(s)\}$, $L_2=\{y(t)\}$. We have
\begin{equation}\label{eq:F_L-linking}
 F_L:S^1\times S^1\xrightarrow{\ L_1\times L_2\ } \Cnf_2(\R^3) \stackrel{r}{\longrightarrow} S^2,\qquad (L_1\times L_2)(s,t)=(x(s), y(t))\ .
\end{equation}
where $r(x,y)=\frac{x-y}{\|x-y\|}$ is a retraction of $\Cnf_2(\R^3)$ onto $S^2$. The Gauss linking number formula \cite{Gauss1833} reads
\begin{equation}\label{eq:linking-gauss}
 \bar{\mu}_{12}(L_1,L_2)=\text{deg}(F_L)=\int_{S^1\times S^1} F^\ast_L(\omega),
\end{equation}
where $\omega\in \Omega^2(S^2)$ is an 
area form of $S^2$.  The first identity in Equation \ref{eq:linking-gauss} is a consequence of the diagrammatic definition of $\bar{\mu}_{12}(L_1,L_2)$ (see \cite{Prasolov-Sossinsky97}) equivalent to the intersection theory definition provided in Section \ref{sec:2}. 

In the following we focus exclusively on relation between $\kappa$-invariants and $\bar{\mu}$-invariants in the $3$-component case. In the context of results \cite{Koschorke04, Koschorke91, Koschorke97, Kom-unlinked09, Kom-Milnor08}
consider a 3-component link $L=\{L_1, L_2, L_3\}$ in $S^3$ 
parametrized by $\{x(s),y(t),z(u)\}$ and 
\begin{equation}\label{eq:F}
 \begin{split}
 & \tilde{F}_L:S^1\times S^1\times S^1\xrightarrow{\ L_1\times L_2\times L_3\ } \Cnf_3(S^3)\stackrel{H}{\longrightarrow} S^2,\\
  & (L_1\times L_2\times L_3)(s,t,u)=(x(s), y(t), z(u)),
  \end{split}
\end{equation}
where $H$ is a projection on the second factor of $\Cnf_3(S^3)\cong S^3\times S^2$, and $\cong$ denotes the homotopy equivalence. The map $H$ may be defined with help of the quaternionic structure of $S^3$ as follows  
\begin{gather}\label{eq:hermans-map}
 \Cnf_3(S^3) \ni (x,y,z) \stackrel{H}{\longrightarrow}
 \frac{\text{pr}(x^{-1}\cdot y)-\text{pr}(x^{-1}\cdot z)}{\|\text{pr}(x^{-1}\cdot y)-\text{pr}(x^{-1}\cdot z)\|}\in S^2,
\end{gather}
where $\cdot$ stands for the quaternionic multiplication, $\ ^{-1}$
is the quaternionic inverse, and $\text{pr}:S^3-\{1\}\longrightarrow \R^3$ the
stereographic projection from $1$, cf. \cite{Kom-Milnor08}.
\begin{theorem}[\cite{Kom-Milnor08}, for $(b)$ also see \cite{Kom-unlinked09}]\label{th:milnor-hopf-S3}
 Let $\tilde{L}=\{L_1,L_2,L_3\}$ be a 3-component link in $S^3$, then the 
 associated map to $\tilde{F}_{L}$  defined in \eqref{eq:F} satisfies
 \begin{itemize}
 \item[$(a)$]  $\text{\rm deg}(\tilde{F}_{L}|_{S^1_i\times S^1_j})=\bar{\mu}_{12}(L_i,L_j)$,
 \item[$(b)$]  whenever $\tilde{L}$ is Borromean, $\tilde{F}_L$ is homotopic to $\pm\, 2\bar{\mu}_{123}\times$ 
 the Hopf map, where the sign depends on the orientation of components. 
 \item[$(c)$] in the general case 
 \[
  \nu(\tilde{F}_L)=\pm\, 2\,\bar{\mu}_{123}(L) \mod\  2\,\text{\rm gcd}(\bar{\mu}_{12}(L_1,L_2), \bar{\mu}_{12}(L_2,L_3), \bar{\mu}_{12}(L_1,L_3)),
  \]
  where $\nu(\tilde{F}_L)$ is the Pontryagin invariant of $\tilde{F}_L$ i.e. 
  the framing of the inverse image of a regular value of $\tilde{F}_L$ 
  (consult \cite{Kom-Milnor08} for a precise definition).
 \item[$(d)$] for Borromean $\tilde{L}$ we have the following formula 
 \[
   \bar{\mu}_{123}=\pm\,\frac{1}{2}\int_{(S^1)^3} \tilde{F}^\ast_L\omega\wedge \eta,
 \]
 where $\omega$ is the area form on $S^2$ and $d\eta=\omega$.
 \end{itemize}
\end{theorem}

\no In the next theorem, Theorem \ref{th:milnor-hopf-R3}, $(a)$, $(b)$ and $(d)$ are extended to link maps valued in $\R^3$.  The theorem has been obtained earlier by Koschorke as Corollary 6.2 in \cite[p. 314]{Koschorke97}, which treats the general $n$-component Borromean case. In the appendix of this article we show how Theorem \ref{th:milnor-hopf-R3} follows from Theorem \ref{th:milnor-hopf-S3}. (Note that $(c)$ of Theorem \ref{th:milnor-hopf-S3} is an original contribution of \cite{Kom-Milnor08}.) Paraphrasing Corollary 6.2 of \cite[p. 314]{Koschorke97} we state

\begin{theorem}[\cite{Koschorke97}]\label{th:milnor-hopf-R3}
 Let $L=\{L_1,L_2,L_3\}$ be a 3-component link in $\R^3$, then the 
 map $F_{L}$ defined in \eqref{eq:link-map}
 satisfies 
 \begin{itemize}
 \item[$(i)$]  $\text{\rm deg}(F_{L}|_{S^1_i\times S^1_j})=\bar{\mu}_{12}(L_i,L_j)$.
 \item[$(ii)$] whenever $L$ is Borromean $F_L$ is homotopic to $\bar{\mu}_{123}\times$ one of the Whitehead products 
 \[
  [\alpha_{3,2},\alpha_{3,1}], [\alpha_{3,1},\alpha_{2,1}], -[\alpha_{3,2},\alpha_{2,1}].
 \]
 \end{itemize}
\end{theorem}

\subsection{Proof of Key Lemma.}

Thanks to $(i)$, in the Borromean case $F_L:(S^1)^2\mapsto\Cnf_3(\R^3)$ is homotopic to a map which is constant when restricted to the $2$-skeleton of the $3$-torus $(S^1)^3$. Therefore, we may integrate directly 
over the torus (similar argument is presented in \eqref{eq:I-homotopy-invariance}), obtaining  by Proposition \ref{th:integral-whitehead}:
\begin{equation}\label{eq:I-mu}
\begin{split}
  \bar{\mu}_{123}(L_1,L_2,L_3) & =\pm \mathcal{I}(F_L)\\
  & =\int_{(S^1)^3} \bigl(F^\ast_L\omega_{1,2}\wedge \eta_{2,3} +F^\ast_L\omega_{2,3}\wedge \eta_{3,1}+F^\ast_L\omega_{3,1}\wedge \eta_{1,2}-F^\ast_L\phi_{1,2,3}\bigr),
 \end{split}
 \end{equation}
where $\eta_{i,j}$ satisfy $d\eta_{i,j}=F^\ast_L\omega_{i,j}$. This proves Key Lemma.
\hfill$\Box$

{\tiny
}

\medskip
\section*{Appendix: Proof of Theorem \ref{th:milnor-hopf-R3}}

\begin{proof}
 The proof of $(i)$ is immediate from the Gauss formula for the linking number in \eqref{eq:linking-gauss}. 
 To prove $(ii)$ consider the  inclusion 
 \[
  j:\R^3\longrightarrow \R^3\cup \{\infty\}\cong S^3,
 \]
 defined by the inverse of stereographic projection and observe that $j$ leads to the inclusion on configuration spaces
 \[
  \widehat{j}:\Cnf_3(\R^3)\longrightarrow \Cnf_3(S^3)\ .
 \]
\no As a first step, we calculate the induced homomorphism cf. \cite[Theorem 2.2, p.34]{Fadell-Husseini01}
\[
 \pi_2(\widehat{j}):\pi_2(\Cnf_3(\R^3))\longrightarrow\pi_2(\Cnf_3(S^3)).
\]
Recall that $\pi_2(\Cnf_3(\R^3))\cong \Z\oplus\Z\oplus\Z$ is generated by classes $\alpha_{i,j}$ represented by $A_{i,j}$, Equation \ref{eq:Aij-def}, and $\pi_2(\Cnf_3(S^3))\cong \Z$ is generated just by one class $\alpha$ represented by $S^2$ factor in $\Cnf_3(S^3)\cong S^3\times S^2$. Let us denote 
\begin{equation}\label{eq:beta-ij}
 \beta_{i,j}:=\pi_2(\widehat{j})(\alpha_{i,j})=a_{i,j}\,\alpha,\qquad a_{i,j}\in \Z,
\end{equation}
 we must determine the coefficients $a_{ij}$.

The symmetric group $\Sigma_3$ acts on 
both $\Cnf_3(\R^3)$ and $\Cnf_3(S^3)$ by permuting the coordinate factors, and the inclusion map $j$ is equivariant with respect to this action. For our purposes we only need to determine the action of the transposition $(2,3)$ on $\alpha_{3,1}$ and $\alpha_{3,2}$.  We claim the following identities:
\begin{equation}\label{eq:alpha_ij-perms}
 (2,3)\alpha_{3,2}=-\alpha_{3,2}\qquad (2,3)\alpha_{3,1}=\alpha_{2,1}\ .
\end{equation}
To justify, we calculate on representatives in \eqref{eq:Aij-def}
\[
\begin{split}
 (2,3) A_{3,2}:\xi\longrightarrow (2,3)(q_1,q_2,q_2+\xi)=(q_1,q_2+\xi, q_2)\\
 (2,3) A_{3,1}:\xi\longrightarrow (2,3)(q_1,q_2,q_1+\xi)=(q_1,q_1+\xi,q_2)\ .
\end{split}
\]
Consider the homotopy
\[
 G:t\longrightarrow (q_1,q_2+(1-t)\xi, q_2-t\,\xi),\qquad t\in [0,1]
\]
Notice that $G$ is well defined in $\Cnf_3(\R^3)$, because $q_2+(1-t)\xi=q_2-t\,\xi$ implies $\xi=0$ 
which contradicts $\xi\in S^2$. This homotopy connects $G_0=(q_1,q_2+\xi, q_2)=(2,3)A_{3,2}$ and 
$G_1=(q_1,q_2, q_2-\xi)=-A_{3,2}$, which shows the first identity in \eqref{eq:alpha_ij-perms}.	
To see the second identity consider the homotopy: $G':t\longrightarrow (q_1,q_1+\xi, (1-t)q_2+t q_3)$, $t\in [0,1]$.

 Recall the fibration Diagram \eqref{eq:fibration-pi}, the image of the fiber $\R^3-\{q_1,q_2\}$ under inclusion $j$ is in $(\R^3\cup \{\infty\})-\{q_1,q_2\}\subset \Cnf_3(S^3)$
implying the following relation
\[
 \beta_{3,1}+\beta_{3,2}=0,\qquad \text{in}\quad \pi_2(\Cnf_3(S^3))\ .
\]
\no Applying \eqref{eq:alpha_ij-perms} to the above equation we obtain $\beta_{2,1}-\beta_{3,2}=0$,
hence
\[
 \beta_{2,1}=\beta_{3,2}=-\beta_{3,1}\qquad \text{in}\quad \pi_2(\Cnf_3(S^3))\ .
\]
\no Thanks to the map defined in \eqref{eq:hermans-map} we observe $\beta_{2,1}=\alpha$ and
\[
 \pi_2(\widehat{j})(\alpha_{2,1})=\alpha,\quad \pi_2(\widehat{j})(\alpha_{3,2})=\alpha,\quad \pi_2(\widehat{j})(\alpha_{3,1})=-\alpha.
\]
Therefore, coefficients in \eqref{eq:beta-ij} are $a_{2,1}=a_{3,2}=-a_{3,1}=1$.
Let $h:S^3\mapsto \Cnf_3(S^3)\cong S^3\times S^2$ be a map such that $p_1\circ h$ is null and $p_2\circ h$ is homotopic to the the Hopf map (where $p_i$ is the projection onto the $i$th factor in $S^3\times S^2$). By naturality of the Whitehead product  \cite[p. 473]{Whitehead78} we obtain from Equation \eqref{eq:2hopf}
\begin{equation}\label{eq:pi-j}
\begin{split}
 \pi_3(\widehat{j})([\alpha_{3,2},\alpha_{3,1}]) & =[\alpha,-\alpha]=-2\,[h],\qquad
 \pi_3(\widehat{j})([\alpha_{3,2},\alpha_{2,1}])=[\alpha,\alpha]=2\,[h],\\
 \pi_3(\widehat{j})([\alpha_{3,1},\alpha_{2,1}]) & =[-\alpha,\alpha]=-2\,[h].
\end{split}
\end{equation}
 In $\pi_3(\Cnf_3(\R^3))$ the following identity, known as Yang-Baxter relation, \cite{Kohno02}, holds
\begin{equation}\label{eq:Y-B}
 [\alpha_{3,2},\alpha_{3,1}+\alpha_{3,2}]=0.
\end{equation}
To justify the identity just consider $\phi:S^2\times S^2\longrightarrow \Cnf_3(\R^3)$, $\phi(\xi_1,\xi_2)=(q_1,q_1+\xi_1,q_1+5\,\xi_2)$, \cite{Fadell-Husseini01}. It is easy to see from Equation \eqref{eq:Aij-def} that
\[
 i_1:=\phi|_{S^2\times\{\ast\}}\cong \alpha_{2,1},\qquad i_2:=\phi|_{\{\ast\}\times S^1}\cong \alpha_{3,1}+\alpha_{3,2}.
\]
Similarly, \eqref{eq:Y-B} immediately follows from the naturality of the Whitehead product and the fact that $[i_1,i_2]=0$ in $\pi_3(S^2\times S^2)$ which  is a direct consequence of the definition of the Whitehead product as an attaching map of the $4$-cell to the 2-skeleton $S^2\vee S^2$ in $S^2\times S^2$. Equations \eqref{eq:Y-B} and \eqref{eq:alpha_ij-perms} yield
\begin{equation}\label{eq:whitehead-identities}
 [\alpha_{3,2},\alpha_{3,1}]=[\alpha_{3,1},\alpha_{2,1}]=-[\alpha_{3,2},\alpha_{2,1}]\ .
\end{equation}
\no Next, we need to prove that the link map $F_L$ associated to a Borromean link $L$ is a multiple of the above Whitehead products in $\pi_3(\Cnf_3(\R^3))$. Because $\bar{\mu}_{12}(L_i,L_j)=0$, we may homotopy $F_L$  to a map $\tilde{f}_L$ constant on the $2$-skeleton of $(S^1)^3$. The map $\tilde{f}_L$ represents an element in $\pi_3(\Cnf_3(\R^3))$, such that
\begin{equation}\label{eq:F_L-in-ker}
 [\tilde{f}_L]\in \ker(\pi_3(\Pi)),\qquad \pi_3(\Pi):\pi_3(\Cnf_3(\R^3))\longrightarrow \pi_3((\Cnf_2(\R^3))^3),
\end{equation}
where $\Pi=\Pi_1\times \Pi_2\times \Pi_3$, and $\Pi_i$ were defined in \eqref{eq:Pi_i}. Recall that 
$\pi_3(\Cnf_3(\R^3))\cong \pi_3(S_{2,1})\oplus\pi_3(S_{3,1}\vee S_{3,2})$, \cite[p. 189]{Whitehead78}, and the group $\pi_3(S_{2,1})\cong \Z$ 
is generated by a Hopf map $h_{21}$. The group $\pi_3(S_{3,1}\vee S_{3,2})\cong \Z\oplus \Z\oplus \Z$
is generated by Hopf maps $h_{3,1}$, $h_{3,2}$ and the Whitehead product $[\alpha_{3,2},\alpha_{3,1}]$
(which immediately follows from the split short exact sequence $0\mapsto \pi_{4}(S^2\times S^2, S^2\vee S^2)\mapsto \pi_3(S^2\vee S^2)\mapsto \pi_3(S^2\times S^2)\mapsto 0$, \cite[p. 492]{Whitehead78}). 
We expand
\[ 
 [\tilde{f}_L]=c_1 [h_{2,1}]+c_2 [h_{3,1}]+c_3 [h_{3,2}]+c_4 [\alpha_{3,2},\alpha_{3,1}],\qquad c_k\in \Z\ .
\]
Since $\pi_3(\Pi_i)([h_{i,j}])=[h_{i,j}]$ for all $1\leq j< i\leq 3$, \eqref{eq:F_L-in-ker} tells us 
that $c_1=c_2=c_3=0$ and 
\[ 
 \tilde{f}_L=c_4 [\alpha_{3,2},\alpha_{3,1}].
\]
\no From Theorem \ref{th:milnor-hopf-S3} and Identities \eqref{eq:pi-j} we conclude
\[
 \pm 2 \bar{\mu}_{123} [h]=\pi_3(\widehat{j})(\tilde{f}_L)= 2\,c_4\,[h],\quad\text{and}\quad  c_4 = \pm\bar{\mu}_{123}\ .
\]
The claim $(ii)$ follows. 

\end{proof}

\begin{thebibliography}{10}

\bibitem{Akhmetiev05}
P.~Akhmetiev.
\newblock On a new integral formula for an invariant of 3-component oriented
  links.
\newblock {\em J. Geom. Phys.}, 53(2):180--196, 2005.

\bibitem{Arnold86}
V.~Arnold.
\newblock The asymptotic {H}opf invariant and its applications.
\newblock {\em Selecta Math. Soviet.}, 5(4):327--345, 1986.
\newblock Selected translations.

\bibitem{Khesin98}
V.~Arnold and B.~Khesin.
\newblock {\em Topological methods in hydrodynamics}, volume 125 of {\em
  Applied Mathematical Sciences}.
\newblock Springer-Verlag, New York, 1998.

\bibitem{Baader07}
S.~Baader.
\newblock Asymptotic {R}asmussen invariant.
\newblock {\em C. R. Math. Acad. Sci. Paris}, 345(4):225--228, 2007.

\bibitem{Baader-Marche07}
S.~Baader and J.~Marche.
\newblock Asymptotic vassiliev invariants for vector fields.
\newblock {\em arXiv:0810.3870}, 2008.

\bibitem{Becker81}
M.~E. Becker.
\newblock Multiparameter groups of measure-preserving transformations: a simple
  proof of {W}iener's ergodic theorem.
\newblock {\em Ann. Probab.}, 9(3):504--509, 1981.

\bibitem{Berger90}
M.~Berger.
\newblock Third-order link integrals.
\newblock {\em J. Phys. A}, 23(13):2787--2793, 1990.

\bibitem{Cantarella00}
J.~Cantarella.
\newblock A general mutual helicity formula.
\newblock {\em R. Soc. Lond. Proc. Ser. A Math. Phys. Eng. Sci.},
  456(2003):2771--2779, 2000.

\bibitem{Kom-hierarchy}
J.~Cantarella, R.~Komendarczyk, and J.~Parsley.
\newblock Higher helicities, rope length and energy.
\newblock {\em in preparation}.

\bibitem{Moffatt95}
A.~Y.~K. Chui and H.~K. Moffatt.
\newblock The energy and helicity of knotted magnetic flux tubes.
\newblock {\em Proc. Roy. Soc. London Ser. A}, 451(1943):609--629, 1995.

\bibitem{Cohen-Lada-May76}
F.~Cohen, T.~J. Lada, and J.~P. May.
\newblock {\em The homology of iterated loop spaces}.
\newblock Springer-Verlag, Berlin, 1976.
\newblock Lecture Notes in Mathematics, Vol. 533.

\bibitem{Kom-Milnor08}
D.~DeTurck, H.~Gluck, R.~Komendarczyk, P.~Melvin, C.~Shonkwiler, and D.~S.
  Vela-Vick.
\newblock Triple linking numbers, ambiguous {H}opf invariants and integral
  formulas for three-component links.
\newblock {\em Mat. Contemp.}, 34:251--283, 2008.

\bibitem{Fadell-Husseini01}
E.~R. Fadell and S.~Y. Husseini.
\newblock {\em Geometry and topology of configuration spaces}.
\newblock Springer Monographs in Mathematics. Springer-Verlag, Berlin, 2001.

\bibitem{Freedman91-2}
M.~Freedman and Z.~He.
\newblock Divergence-free fields: energy and asymptotic crossing number.
\newblock {\em Ann. of Math. (2)}, 134(1):189--229, 1991.

\bibitem{Gambaudo-Ghys97}
J.-M. Gambaudo and {\'E}.~Ghys.
\newblock Enlacements asymptotiques.
\newblock {\em Topology}, 36(6):1355--1379, 1997.

\bibitem{Gauss1833}
C.~F. Gauss.
\newblock Integral formula for linking number.
\newblock {\em Zur Mathematischen Theorie der Electrodynamische Wirkungen
  ({C}ollected {W}orks, Vol. 5), Koniglichen Gesellschaft des Wissenschaften,
  G{\"o}ttingen}, 2ed.(3):605, 1833.

\bibitem{Haefliger78}
A.~Haefliger.
\newblock Whitehead products and differential forms.
\newblock In {\em Differential topology, foliations and {G}elfand-{F}uks
  cohomology ({P}roc. {S}ympos., {P}ontif\'\i cia {U}niv. {C}at\'olica, {R}io
  de {J}aneiro, 1976)}, volume 652 of {\em Lecture Notes in Math.}, pages
  13--24. Springer, Berlin, 1978.

\bibitem{Hain84}
R.~Hain.
\newblock Iterated integrals and homotopy periods.
\newblock {\em Mem. Amer. Math. Soc.}, 47(291):iv+98, 1984.

\bibitem{Khesin05}
B.~Khesin.
\newblock Topological fluid dynamics.
\newblock {\em Notices Amer. Math. Soc.}, 52(1):9--19, 2005.

\bibitem{Kohno02}
T.~Kohno.
\newblock Loop spaces of configuration spaces and finite type invariants.
\newblock In {\em Invariants of knots and 3-manifolds (Kyoto, 2001)}, volume~4
  of {\em Geom. Topol. Monogr.}, pages 143--160 (electronic). Geom. Topol.
  Publ., Coventry, 2002.

\bibitem{Kom-unlinked09}
R.~Komendarczyk.
\newblock The third order helicity of magnetic fields via link maps.
\newblock {\em Comm. Math. Phys.}, 292(2):431, 2009.

\bibitem{Koschorke91}
U.~Koschorke.
\newblock Link homotopy with many components.
\newblock {\em Topology}, 30(2):267--281, 1991.

\bibitem{Koschorke97}
U.~Koschorke.
\newblock A generalization of {M}ilnor's {$\mu$}-invariants to
  higher-dimensional link maps.
\newblock {\em Topology}, 36(2):301--324, 1997.

\bibitem{Koschorke04}
U.~Koschorke.
\newblock Link homotopy in {$S\sp n\times \mathbb R\sp {m-n}$} and higher order
  {$\mu$}-invariants.
\newblock {\em J. Knot Theory Ramifications}, 13(7):917--938, 2004.

\bibitem{Laurence-Stredulinsky00b}
P.~Laurence and E.~Stredulinsky.
\newblock Asymptotic {M}assey products, induced currents and {B}orromean torus
  links.
\newblock {\em J. Math. Phys.}, 41(5):3170--3191, 2000.

\bibitem{Laurence-Stredulinsky00a}
P.~Laurence and E.~Stredulinsky.
\newblock A lower bound for the energy of magnetic fields supported in linked
  tori.
\newblock {\em C. R. Acad. Sci. Paris S\'er. I Math.}, 331(3):201--206, 2000.

\bibitem{Mayer03}
C.~Mayer.
\newblock Topological link invariants of magnetic fields.
\newblock {\em Ph.D. thesis}, 2003.

\bibitem{Mellor03}
B.~Mellor and P.~Melvin.
\newblock A geometric interpretation of {M}ilnor's triple linking numbers.
\newblock {\em Algebr. Geom. Topol.}, 3:557--568 (electronic), 2003.

\bibitem{Milnor54}
J.~Milnor.
\newblock Link groups.
\newblock {\em Ann. of Math. (2)}, 59:177--195, 1954.

\bibitem{Milnor57}
J.~Milnor.
\newblock Isotopy of links.
\newblock In R.~Fox, editor, {\em Algebraic Geometry and Topology}, pages
  280--306. Princeton University Press, 1957.

\bibitem{Moffatt85}
H.~K. Moffatt.
\newblock Magnetostatic equilibria and analogous {E}uler flows of arbitrarily
  complex topology. {I}. {F}undamentals.
\newblock {\em J. Fluid Mech.}, 159:359--378, 1985.

\bibitem{Novikov88}
S.~P. Novikov.
\newblock Analytical theory of homotopy groups.
\newblock In {\em Topology and geometry---Rohlin Seminar}, volume 1346 of {\em
  Lecture Notes in Math.}, pages 99--112. Springer, Berlin, 1988.

\bibitem{Porter80}
R.~Porter.
\newblock Milnor's {$\bar \mu $}-invariants and {M}assey products.
\newblock {\em Trans. Amer. Math. Soc.}, 257(1):39--71, 1980.

\bibitem{Prasolov-Sossinsky97}
V.~V. Prasolov and A.~B. Sossinsky.
\newblock {\em Knots, links, braids and 3-manifolds}, volume 154 of {\em
  Translations of Mathematical Monographs}.
\newblock American Mathematical Society, Providence, RI, 1997.
\newblock An introduction to the new invariants in low-dimensional topology,
  Translated from the Russian manuscript by Sossinsky [Sosinski{\u\i}].

\bibitem{Priest84}
E.~Priest.
\newblock {\em Solar Magnetohydrodynamics}.
\newblock D.Redidel Publishing Comp., 1984.

\bibitem{Schwarz95}
G.~Schwarz.
\newblock {\em Hodge decomposition---a method for solving boundary value
  problems}, volume 1607 of {\em Lecture Notes in Mathematics}.
\newblock Springer-Verlag, Berlin, 1995.

\bibitem{Sinha-Walter08}
D.~Sinha and B.~Walter.
\newblock Lie coalgebras and rational homotopy theory {II}: Hopf invariants.
\newblock {\em arXiv.org:0809.5084}, 2008.

\bibitem{Sullivan77}
D.~Sullivan.
\newblock Infinitesimal computations in topology.
\newblock {\em Inst. Hautes \'Etudes Sci. Publ. Math.}, (47):269--331 (1978),
  1977.

\bibitem{Hornig04}
H.~v.~Bodecker and G.~Hornig.
\newblock Link invariants of electromagnetic fields.
\newblock {\em Phys. Rev. Lett.}, 92(3):030406, 4, 2004.

\bibitem{Verjovsky94}
A.~Verjovsky and R.~F. Vila~Freyer.
\newblock The {J}ones-{W}itten invariant for flows on a {$3$}-dimensional
  manifold.
\newblock {\em Comm. Math. Phys.}, 163(1):73--88, 1994.

\bibitem{Vogel03}
T.~Vogel.
\newblock On the asymptotic linking number.
\newblock {\em Proc. Amer. Math. Soc.}, 131(7):2289--2297 (electronic), 2003.

\bibitem{Whitehead78}
G.~W. Whitehead.
\newblock {\em Elements of homotopy theory}, volume~61 of {\em Graduate Texts
  in Mathematics}.
\newblock Springer-Verlag, New York, 1978.

\bibitem{Woltjer58}
L.~Woltjer.
\newblock A theorem on force-free magnetic fields.
\newblock {\em Proc. Nat. Acad. Sci. U.S.A.}, 44:489--491, 1958.

\end{thebibliography}
\end{document}